\documentclass[english,aip,cha]{revtex4-1}
\usepackage[T1]{fontenc}
\usepackage[latin9]{inputenc}
\usepackage{amsthm}
\usepackage{amsmath}
\usepackage{amssymb}
\usepackage{graphicx}
\usepackage{natbib}

\usepackage{babel}
\makeatletter
\@ifundefined{textcolor}{}
{%
 \definecolor{BLACK}{gray}{0}
 \definecolor{WHITE}{gray}{1}
 \definecolor{RED}{rgb}{1,0,0}
 \definecolor{GREEN}{rgb}{0,1,0}
 \definecolor{BLUE}{rgb}{0,0,1}
 \definecolor{CYAN}{cmyk}{1,0,0,0}
 \definecolor{MAGENTA}{cmyk}{0,1,0,0}
 \definecolor{YELLOW}{cmyk}{0,0,1,0}
}
\numberwithin{equation}{section}
\numberwithin{figure}{section}
\theoremstyle{plain}
\newtheorem{thm}{\protect\theoremname}
  \theoremstyle{plain}
  \newtheorem{lem}[thm]{\protect\lemmaname}
  \theoremstyle{plain}
  \newtheorem{prop}[thm]{\protect\propositionname}
  \theoremstyle{remark}
  \newtheorem{rem}[thm]{\protect\remarkname}

\makeatother

\usepackage{babel}
  \providecommand{\lemmaname}{Lemma}
  \providecommand{\propositionname}{Proposition}
  \providecommand{\remarkname}{Remark}
\providecommand{\theoremname}{Theorem}

\begin{document}

\title{Amplitude equations and fast transition to chaos in rings of coupled
oscillators}

\author{S. Yanchuk}
\affiliation{ Institute of Mathematics, Humboldt University at Berlin,
Unter den Linden 6, 10099 Berlin, Germany}

\author{P. Perlikowski} 
\affiliation{ Division of Dynamics, Technical University of Lodz, 90-924
Lodz, Poland}

\author{M. Wolfrum}
\affiliation{ Weierstrass Institute for Applied Analysis and Stochastics,
Mohrenstrasse 39, 10117 Berlin, Germany}

\author{A. Stefa\'{n}ski}
\affiliation{ Division of Dynamics, Technical University of Lodz, 90-924
Lodz, Poland}

\author{and T.~Kapitaniak}
\affiliation{ Division of Dynamics, Technical University of Lodz, 90-924
Lodz, Poland}

\begin{abstract}
We study the coupling induced destabilization in an array of identical
oscillators coupled in a ring structure where the number of oscillators
in the ring is large. The coupling structure includes different types
of interactions with several next neighbors. We derive an amplitude
equation of Ginzburg-Landau type, which describes the destabilization
of a uniform stationary state in a ring with a large number of nodes.
Applying these results to unidirectionally coupled Duffing oscillators,
we explain the phenomenon of a fast transition to chaos, which has
been numerically observed in such systems. More specifically, the
transition to chaos occurs on an interval of a generic control parameter
that scales as the inverse square of the size of the ring, i.e. for
sufficiently large system, we observe practically an immediate transition
to chaos.
\end{abstract}


\maketitle
34C15, 34C23, 34C28, 37G35

\begin{quotation}
The understanding of the dynamical behavior of networks of
coupled oscillators can contribute to the explanation of various collective
phenomena that can be observed in coupled systems in biology, economy,
or physics \cite{Strogatz2001,Mosekilde2002,Pikovsky2001}. Discrete
media, e.g. in neural systems, can exhibit complex coupling structures
including feed-forward loops, large coupling ranges, and interaction
mechanisms of different kind. In such systems, the influence of the
coupling can transform the simple dynamics of a single oscillator
into complicated spatio-temporal structures of the network dynamics.
Amplitude equations of Ginzburg-Landau type have been a powerful tool
for a universal description of spatio-temporal dynamics and pattern
formation in continuous media. In this paper, we apply this technique
to a large class of coupled oscillator systems with a ring structure,
where we assume that the number of oscillators in the ring is large.
In addition, we use the amplitude equation together with the corresponding
scaling laws to explain the coupling induced fast transition from
homogeneous stationary behavior to high dimensional chaos, which can
be observed in certain coupled oscillator systems.
\end{quotation}

\section{Introduction }

Networks of coupled oscillators have been the subject of extensive
research in the last decade. Coupled systems can display a huge variety
of dynamical phenomena, starting from synchronization phenomena in
various types of inhomogeneous or irregular networks, up to complex
collective behavior, such as for example various forms of phase transitions,
traveling waves \cite{Kumar2010,Ermentrout2001,Popovych2011}, phase-locked
patterns, amplitude death states \cite{Dodla2004}, or so called chimera
states that display a regular pattern of coherent and incoherent motion
\cite{Abrams2004,Omel2008,Sethia2008,Wolfrum2011}. Of particular
interest are situations, where complex spatio-temporal structures
can emerge in regular arrays of identical units induced only by the
coupling interaction. In many cases, the resulting phenomena differ
substantially from corresponding situations in continuous media \cite{Nakao2010}
and depend strongly on the underlying network topology.

Our specific interest is in the emergence of spatio-temporal structures
in a homogeneous array of identical units that have a stable uniform
equilibrium at which the coupling vanishes. As a classical paradigm,
the Turing instability gives an example of a coupling induced instability
in such a setting. This phenomenon has of course a direct counterpart
in the discrete setting, but it turns out that there appear some genuinely
new phenomena. In \cite{Yanchuk2008a,Perlikowski2010,Yanchuk2011}
it has been shown that also in a ring of unidirectionally coupled
oscillators, i.e. in a purely convective setting, the Eckhaus scenario
of coexisting diffusive patterns can be observed. In \cite{Perlikowski2010a}
it has been shown that Duffing oscillators coupled in the same way,
exhibit a complex transition to spatio-temporal chaos. In this paper
we develop a general theoretical framework for such phenomena in large
arrays. We derive an amplitude equation of Ginzburg-Landau type that
governs the local dynamics close to the destabilization of the uniform
steady state. It resembles several features that are already well
known in the context of reaction-diffusion systems \cite{Eckhaus1965,Tuckerman1991}.
But in contrast to these results, it can be applied to much more general
coupling mechanisms, including also the case of unidirectional and
anti-diffusive interaction and allows also for a mixture of such interactions
with several next neighbors. Such an interplay of attractive and repulsive
coupling is an essential feature of neural systems. As a specific
feature, the convective part will appear in the amplitude equation
as a rotation of the coordinates in an intermediate time scale that
is faster than the diffusive processes described by the Ginzburg-Landau
equation.

Having deduced the amplitude equation and the corresponding scaling
laws in terms of the number of oscillators, which is assumed to tend
to infinity, we use this theory for the explanation of a specific
phenomenon that has been reported in \cite{Perlikowski2010a} for
a ring of unidirectionally coupled Duffing oscillators. There, it
has been shown numerically that for a large number of oscillators
$N$, there is an almost immediate transition from homogeneous stationary
behavior to high-dimensional chaos. Based on our amplitude equation,
we argue that in such systems, one can expect generically that such
a transition occurs within a parameter interval of the size $1/N^{2}$.
We consider a generic case, where the control parameter enters already
the linear parts of the dynamical equations, e.g. a diffusive coupling
strength. Finally, we demonstrate this phenomenon by a numerical example,
where we also evaluate the scaling behavior of the parameter interval
where the transition to chaos takes place for an increasing number
of oscillators.

\section{Model equation, spectral conditions, and notations }

\begin{figure}[h]
\begin{centering}
\includegraphics[scale=0.75]{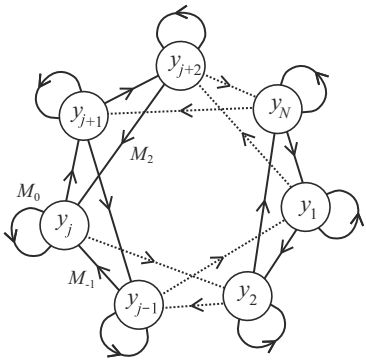} 
\par\end{centering}

\caption{An example of a ring of $N$ coupled oscillators. Apart from the self-coupling
$M_{0},$ each oscillator $y_{j}$ is also coupled with $y_{j+2}$
($M_{2}$) as well as $y_{j-1}$ ($M_{-1}$). See Eq. (\ref{eq:main})
for the equation of motion. \label{fig:ring}}
\end{figure}
We are interested in a system of $N$ identical coupled oscillators
that has an uniform equilibrium, where the coupling vanishes. The
coupling network is organized in a ring structure, where interactions
of several next neighbors are possible. Such systems can be written
in general form as 
\begin{equation}
\dot{y}_{j}=\sum_{m}M_{m}(p)y_{j+m}+h(y_{j},y_{j+1},\dots y_{j+N};p),\label{eq:main}
\end{equation}
where $y_{j}\in\mathbb{R}^{n},\; j=1,\dots,N$ describes the state
of the $j$-th oscillator. The ring structure is induced by periodic
boundary conditions, i.e. all indexes have to be considered modulo
$N$. The linear part of the dynamics is given by the $n\times n$
matrices $M_{m}(p),\; m=1,\dots,N$, depending on the bifurcation
parameter $p$, which account for the coupling to the m-th neighbor;
in particular $M_{0}(p)$ describes the linear part of the local dynamics
(self-coupling). The nonlinear part $h$, again including a local
dependence and a dependence on the $m$-th neighbor, should vanish
at the origin $h(0,\dots,0;p)=0$ and have also zero derivatives there.
Note that this system is symmetric (equivariant) with respect to index
shift. Figure \ref{fig:ring} illustrates an example with self coupling
and coupling to the neighbor on the left and to the second neighbor
on the right. The specific form of (\ref{eq:main}) also implies that
the coupling vanishes at the equilibrium $y_{1}=\cdots=y_{N}=0$,
which is true e.g. when the coupling is a function of the difference
$y_{j}-y_{m}$ for any two coupled oscillators $j,m$.

The characteristic equation for the linearization at the zero equilibrium
of (\ref{eq:main}) can be factorized as 
\[
\chi(p,\lambda,e^{i2\pi j/N})=\det\left[\lambda\mathrm{Id}-\sum_{m}e^{im2\pi j/N}M_{m}(p)\right]=0,
\]
where $\mathrm{Id}$ denotes the identity matrix in $\mathbb{R}^{n}$
and the index $j=1,2,\dots,N$ accounts for the $N$-th roots of unity
that appear as the eigenvalues of the circular coupling structure
\cite{Pecora1998}. Following the approach in \cite{Yanchuk2008a,Perlikowski2010},
we replace for large $N$ the discrete numbers $2\pi j/N$ by a continuous
parameter $\varphi$, and obtain the \emph{asymptotic continuous spectrum}
\begin{equation}
\Lambda_{p}=\left\{ \lambda\in\mathbb{C}:\,\chi(p,\lambda,e^{i\varphi})=\det\left[\lambda\mathrm{\mathrm{Id}}-\sum_{m}e^{im\varphi}M_{m}(p)\right]=0,\,\varphi\in[0,2\pi)\right\} ,\label{eq:ACS}
\end{equation}
which contains all eigenvalues and which for large $N$ is covered
densely by the eigenvalues. Since the expression (\ref{eq:ACS}) is
periodic in $\varphi$, the asymptotic continuous spectrum $\Lambda_{p}$
has generically the form of one or several closed curves $\lambda_{p}(\varphi)$
in the complex plane, parametrized by $\varphi$. 

\begin{figure}
\centering{}\includegraphics[width=0.22\linewidth]{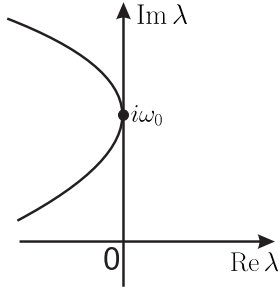}\caption{\label{fig:ACS} Asymptotic continuous spectrum $\Lambda_{p}$ at
the destabilization; schematically.}
\end{figure}

At the bifurcation value $p=0$, we assume that the asymptotic continuous
spectrum touches the imaginary axis at at some point $i\omega_{0}$
(see Fig.~\ref{fig:ACS}), i.e. the following conditions are fulfilled
\begin{equation}
\lambda_{0}(\varphi_{0})=i\omega_{0},\quad\frac{\partial\lambda_{0}}{\partial\varphi}(\varphi_{0})=i\kappa_{1},\quad\kappa_{1}\in\mathbb{R}.\label{eq:bifcond}
\end{equation}
The first condition from (\ref{eq:bifcond}) means that the point
$i\omega_{0}$ belongs to the asymptotic continuous spectrum, while
the second condition from (\ref{eq:bifcond}) guarantees that the
real part $\Re(\lambda(\varphi))$ is tangent to the zero axis at
$\varphi=\varphi_{0}$ (see Fig.~\ref{fig:ACS}). This tangency condition
is quite natural and describes the condition for the destabilization
(or bifurcation) of the zero solution a large ring of coupled oscillators.
Indeed, if the spectrum $\Lambda_{p}$ is contained in the left half
of the complex plane with $\mathrm{Re}\lambda<0$, then the uniform
equilibrium $y_{1}=\cdots=y_{N}=0$ is asymptotically stable. As soon
as $\Lambda_{p}$ crosses the imaginary axis, it becomes unstable
for sufficiently large $N$. 

Before we present our main result, we now introduce some useful notations
and formulate a technical lemma, which follows from the bifurcation
conditions (\ref{eq:bifcond}). With $v_{0}$ and $v_{1}$ we denote
the eigenvector and the adjoint eigenvector to the critical eigenvalue
$\lambda_{0}(\varphi_{0})=i\omega_{0}$, which we assumed in (\ref{eq:bifcond})
to exist for $p=0$. Moreover, it will be convenient to denote by
\[
L_{0}=\sum_{m}e^{im\varphi_{0}}M_{m}(0),
\]
\[
L_{1}=\sum_{m}me^{im\varphi_{0}}M_{m}(0),
\]
\[
L_{2}=\sum_{m}m^{2}e^{im\varphi_{0}}M_{m}(0)
\]
the ``moments'' of the coupling matrices $M_{m}$. In this notation,
the equations for the eigenvectors $v_{0}$ and $v_{1}$ read as 
\begin{equation}
\left[i\omega_{0}\mathrm{Id}-L_{0}\right]v_{0}=0,\label{eq:v0}
\end{equation}
 
\begin{equation}
\left[-i\omega_{0}\mathrm{Id}-L_{0}^{*}\right]v_{1}=0.\label{eq:v1}
\end{equation}
 These vectors can be normalized as 
\begin{equation}
\left|v_{0}\right|^{2}=1,\quad\left\langle v_{0},v_{1}\right\rangle =1.\label{eq:v0v1}
\end{equation}
Finally, we expand the matrices $M_{m}(p)$ with respect to the parameter
$p$ and write them them as $M_{m}(0)+pK_{m}+\mathcal{O}\left(p^{2}\right)$.
In this way, we can define 
\[
L_{K}=\sum e^{im\varphi_{0}}K_{m}.
\]

\begin{lem}
Assume that $\varphi_{0}$ is a regular point of the asymptotic continuous
spectrum (\ref{eq:ACS}), such that $\lambda_{0}(\varphi)$ exists
and is locally differentiable in a small neighborhood of $\varphi_{0}$.
Further, let the bifurcation condition (\ref{eq:bifcond}) hold. Then
\begin{equation}
\left\langle L_{1}v_{0},v_{1}\right\rangle =\kappa_{1},\label{eq:k1}
\end{equation}

\end{lem}
The proof of this Lemma will be given together with the proof of the
main result in Section \ref{sec:Proofs}.

\section{Main Result: Reduction to Ginzburg-Landau Equation}

In this section, we present an amplitude equation that describes for
system (\ref{eq:main}) the dynamics close to the destabilization
threshold in the limit of large $N$. We perform a limiting procedure,
where we assume that a fixed number $R$ of neighbors is involved
in the coupling while the total number of oscillators $N$ tends to
infinity. Also the nonlinearity $h$ is assumed to depend on $R$
next neighbors only. Hence, in the sequel all summation will extend
over $m\in\{-R,\dots,R\}$ and the nonlinearity will be written as
$h(y_{j-R},\dots,y_{j+R};p)$, independently on $N$. In this way,
the coupling becomes local in the limit $N\to\infty$, and coupling
terms will be approximated by derivatives of the amplitude. Consequently,
our results will be valid for large rings, where the coupling range
is small compared to the total size. It will be shown that the amplitude
equation has the form of a complex Ginzburg-Landau equation with periodic
boundary conditions. The derivation is quite technical and we give
here in Proposition 2 only the main assertions. The proof is deferred
to Section \ref{sec:Proofs}.
\begin{prop}
\label{prop2}Assume that the bifurcation condition (\ref{eq:bifcond})
holds and $\varphi_{0}$ is a regular point of the asymptotic continuous
spectrum $\Lambda_{0}$. Additionally, let the points $\pm3i\omega_{0}$
not belong to the asymptotic continuous spectrum with $\varphi=\varphi_{0}$
(nonresonance condition). Let the nonlinearity $h$ be of third order.
Introduce the small parameter
\[
\varepsilon=\frac{1}{N}
\]
and apply the multiple scale ansatz 
\begin{equation}
y_{j}(t)=\varepsilon e^{i\omega_{0}t+i\varphi_{0}j}v_{0}A(T_{1},x_{1},T_{2})+\varepsilon^{3}e^{3i\left(\omega_{0}t+\varphi_{0}j\right)}v_{2}A^{3}(T_{1},x_{1},T_{2})+c.c.,\label{eq:multscale}
\end{equation}
with the amplitude $A\in\mathbb{C}$ depending on the rescaled coordinates
$T_{1}=\varepsilon t,$ $T_{2}=\varepsilon^{2}t,$ and $x_{1}=\varepsilon j$
($c.c.$ denotes complex conjugated terms, $\nu_{2}\in\mathbb{C}^{n}$)
to the following system

\begin{equation}
\dot{y}_{j}=\sum_{m=-R}^{R}\left(M_{m}(0)+\varepsilon^{2}rK_{m}+\mathcal{O}\left(\varepsilon^{4}\right)\right)y_{j+m}+h(y_{j-R},\dots y_{j+R};p),\label{eq:perteq}
\end{equation}
with the rescaled parameter $p=\varepsilon^{2}r$ and $j=1,\dots,N$
with periodic boundary conditions. Then, the solvability conditions
up to the order $\varepsilon^{3}$ imply the following partial differential
equation of Ginzburg-Landau type 
\begin{equation}
\partial_{T_{2}}u=r\kappa_{2}u+\frac{\kappa_{3}}{2}\partial_{\xi}^{2}u+\zeta u\left|u\right|^{2}\label{eq:GL}
\end{equation}
with periodic boundary conditions 
\[
u(\xi,T_{2})=u(\xi+1,T_{2}),
\]
where $u(\xi,T_{2})$ with $\xi\in[0,1]$ is related to the amplitude
$A$ by 
\begin{equation}
A(T_{1},x_{1},T_{2})=u(\kappa_{1}T_{1}+x_{1},T_{2}).\label{eq:Au}
\end{equation}
The coefficient $\kappa_{1}$is given by (\ref{eq:bifcond}) or (\ref{eq:k1}),
and 
\[
\kappa_{2}=\left\langle L_{K}v_{0},v_{1}\right\rangle ,\quad\kappa_{3}=\left\langle L_{2}v_{0},v_{1}\right\rangle .
\]
 Finally, the coefficient of the nonlinearity $\zeta$ and the vector
$v_{2}\in\mathbb{C}^{n}$ have to be determined by the nonlinearity
$h$ according to (\ref{eq:zeta}) and (\ref{eq:v2}).
\end{prop}
According to this proposition, small solutions of a coupled system
of the form (\ref{eq:perteq}) that has a parameter close to criticality,
i.e. $p=O(\varepsilon^{2})$ can be approximated in the form (\ref{eq:multscale}),
where the amplitude $A$ is related to a solution of the Ginzburg-Landau
equation (\ref{eq:GL}) via (\ref{eq:Au}). Note that the relation
(\ref{eq:Au}) introduces a rotating frame on the timescale $T_{1}=\varepsilon t$
with rotation velocity $\kappa_{1}.$ The time evolution of the Ginzburg-Landau
equation enters only on the slowest time scale $T_{2}=\varepsilon^{2}t$. 
\begin{rem}
For the case of a symmetric coupling ($M_{k}=M_{-k}$) the amplitude
equation (\ref{eq:GL}) has all real coefficients $\kappa_{1}$, $\kappa_{2}$
and $\kappa_{3}$ and the rotation velocity vanishes $\kappa_{1}=0$.
\end{rem}

\section{Proofs\label{sec:Proofs}}

\noindent \textbf{Proof of Lemma 1. }We have to show that $\kappa_{1}$,
defined by (\ref{eq:bifcond}), can be calculated as given in (\ref{eq:k1}).
To this end, we differentiate the eigenvalue equation (cf. (\ref{eq:v0}))
\[
\left[\lambda_{0}\left(\varphi\right)\mathrm{Id}-L_{0}\right]v\left(\varphi\right)=0,
\]
with respect to $\varphi$: 
\[
\left[\mathrm{Id}\frac{\partial}{\partial\varphi}\lambda_{0}(\varphi)-i\sum_{m}me^{im\varphi}M_{m}(0)\right]v(\varphi)+\left[\lambda_{0}(\varphi)\mathrm{Id}-\sum_{m}e^{im\varphi}M_{m}(0)\right]\frac{\partial}{\partial\varphi}v(\varphi)=0.
\]
 Evaluating the obtained expression at $\varphi=\varphi_{0},$ $\lambda(\varphi_{0})=i\omega_{0}$,
and $v(\varphi_{0})=v_{0}$, we obtain 
\[
\left[\mathrm{Id}i\kappa_{1}-iL_{1}\right]v_{0}+\left[i\omega_{0}\mathrm{Id}-L_{0}\right]\frac{\partial v}{\partial\varphi}(\varphi_{0})=0.
\]
 The projection onto $v_{0}$ gives 
\[
i\kappa_{1}\left\langle v_{0},v_{1}\right\rangle -i\left\langle L_{1}v_{0},v_{1}\right\rangle +\left\langle \left[i\omega_{0}\mathrm{Id}-L_{0}\right]\frac{\partial v}{\partial\varphi}(\varphi_{0}),v_{1}\right\rangle =0,
\]
 which implies 
\[
i\kappa_{1}-i\left\langle L_{1}v_{0},v_{1}\right\rangle +\left\langle \frac{\partial v}{\partial\varphi}(\varphi_{0})v_{0},\left[-i\omega_{0}\mathrm{Id}-L_{0}^{*}\right]v_{1}\right\rangle =0.
\]
 Taking into account (\ref{eq:v1}), we obtain the relation (\ref{eq:k1}),
which proves the Lemma.

\noindent \textbf{Proof of Proposition 2.} Substituting the multiple
scale ansatz (\ref{eq:multscale}) into (\ref{eq:perteq}), we obtain
\[
\varepsilon\frac{d}{dt}\left(e^{i\omega_{0}t+i\varphi_{0}j}v_{0}A+\varepsilon^{2}e^{3i\left(\omega_{0}t+\varphi_{0}j\right)}v_{2}A^{3}+c.c.\right)
\]
 
\[
=\varepsilon\sum_{m}\left(M_{m}(0)+\varepsilon^{2}rK_{m}\right)e^{i\omega_{0}t+i\varphi_{0}j}e^{im\varphi_{0}}v_{0}A(T_{1},x_{1}+m\varepsilon,T_{2})
\]
\[
+\varepsilon^{3}\sum_{m}\left(M_{m}(0)+\varepsilon^{2}rK_{m}\right)e^{3i\left(\omega_{0}t+\varphi_{0}j\right)}e^{3im\varphi_{0}}v_{2}A^{3}(T_{1},x_{1}+m\varepsilon,T_{2})+c.c.
\]
\[
+h(y_{j-R},\dots,y_{j+R};\varepsilon^{2}r)
\]
Dividing the obtained equation by $\varepsilon$, expanding necessary
arguments of $A$, we obtain up to the terms of the order $\varepsilon^{2}$
(the complex conjugated terms are omitted here for brevity)

\[
e^{i\omega_{0}t+i\varphi_{0}j}v_{0}\left(i\omega_{0}A+\varepsilon\partial_{T_{1}}A+\varepsilon^{2}\partial_{T_{2}}A\right)+3i\omega_{0}\varepsilon^{2}e^{3i\left(\omega_{0}t+\varphi_{0}j\right)}v_{2}A^{3}
\]
\[
=e^{i\omega_{0}t+i\varphi_{0}j}\sum_{m}\left(M_{m}(0)+\varepsilon^{2}rK_{m}\right)e^{im\varphi_{0}}v_{0}\left(A+m\varepsilon\partial_{x_{1}}A+\frac{1}{2}m^{2}\varepsilon^{2}\partial_{x_{1}}^{2}A\right)
\]
\[
+\varepsilon^{2}\sum_{m}M_{m}(0)e^{3i\left(\omega_{0}t+\varphi_{0}j\right)}e^{3im\varphi_{0}}v_{2}A^{3}
\]
\[
+\varepsilon^{2}e^{i\omega_{0}t+i\varphi_{0}j}A\left|A\right|^{2}h_{1}\left(v_{0}\right)+\varepsilon^{2}e^{3i(\omega_{0}t+\varphi_{0}j)}A^{3}h_{2}\left(v_{0}\right)
\]
where $h_{1}\left(v_{0}\right)$ and $h_{2}\left(v_{0}\right)$ are
determined by the leading terms in the expansion of the nonlinearity.
Note that due to our assumption, $h$ is of third order and the leading
order terms are given by expanding $h(y_{j-R},\dots,y_{j+R};0)$ in
the homogeneous state $y_{m}=\alpha v_{0}+c.c,\quad m=j-R,\dots,j+R$
with respect to $\alpha\in\mathbb{C}$. The solvability condition
requires that different harmonics as well as different orders of $\varepsilon$
up to $\varepsilon^{2}$ are equal. We start with the first harmonic.
By equating the terms containing $e^{i\omega_{0}t+i\varphi_{0}j}$
we obtain the following equation

\[
v_{0}\left(i\omega_{0}A+\varepsilon\partial_{T_{1}}A+\varepsilon^{2}\partial_{T_{2}}A\right)=
\]
\[
=\sum_{m}\left(M_{m}(0)+\varepsilon^{2}rK_{m}\right)e^{im\varphi_{0}}v_{0}\left(A+m\varepsilon\partial_{x_{1}}A+\frac{1}{2}m^{2}\varepsilon^{2}\partial_{x_{1}}^{2}A\right)+\varepsilon^{2}A\left|A\right|^{2}h_{1}\left(v_{0}\right).
\]
Since it should be satisfied for all $\varepsilon$, we first consider
the $\varepsilon^{0}$ equation
\[
i\omega_{0}Av_{0}=AL_{0}v_{0},
\]
 which holds according to the spectral condition (\ref{eq:v0}). The
$\varepsilon^{1}$ terms result into
\[
v_{0}\partial_{T_{1}}A-L_{1}v_{0}\partial_{x_{1}}A=0.
\]
Multiplication with $v_{1}^{T}$ and using (\ref{eq:k1}) from Lemma
1, we obtain
\[
\partial_{T_{1}}A-\kappa_{1}\partial_{x_{1}}A=0.
\]
This will be accounted for by introducing the new amplitude $u$ by

\[
u(\xi,T_{2})=u(\kappa_{1}T_{1}+x_{1},T_{2})=A(T_{1},x_{1},T_{2})
\]
in a correspondingly rotating coordinate $\xi=\kappa_{1}T_{1}+x_{1}.$
Finally, the $\varepsilon^{2}$ terms result into
\[
v_{0}\partial_{T_{2}}A=rAL_{K}v_{0}+\frac{1}{2}\partial_{x_{1}}^{2}AL_{2}v_{0}+A\left|A\right|^{2}h_{1}\left(v_{0}\right)
\]
 Note that the dependence on $T_{1}$ does not show up in this equation.
Hence, after multiplication with $v_{1}^{T}$, we can write it in
terms of $u$ as

\[
\partial_{T_{2}}u=r\kappa_{2}u+\frac{\kappa_{3}}{2}\partial_{\xi}^{2}u+\zeta u\left|u\right|^{2},
\]
where 
\begin{equation}
\zeta=\left\langle h_{1}\left(v_{0}\right),v_{1}\right\rangle .\label{eq:zeta}
\end{equation}
 Finally, it is simple to check that the solvability of the terms
for the third harmonic leads to the expression
\begin{equation}
v_{2}=\left[3i\omega_{0}-\sum_{m}M_{m}(0)e^{3im\varphi_{0}}\right]^{-1}h_{2}\left(v_{0}\right)\label{eq:v2}
\end{equation}
Here, the existence of a nonzero solution $v_{2}$ is guaranteed by
the nonresonance condition. Indeed, since the points $\pm3i\omega_{0}$
do not belong to the asymptotic continuous spectrum with $\varphi=\varphi_{0}$,
the matrix $\left[3i\omega_{0}-\sum_{m}M_{m}(0)e^{3im\varphi_{0}}\right]$
is non-singular. Finally, note that the set of equations should be
complemented by periodic boundary conditions in $\xi$, taking into
account the ring structure of system (\ref{eq:perteq}). The proposition
is proved.

\section{Application: scaling of the bifurcation transition}

It is well known that in the complex Ginzburg-Landau equation (\ref{eq:GL})
a destabilization transition is possible, where in a series of bifurcations
a homogeneous stationary state loses its stability and, via periodic
and quasiperiodic motion, a regime of spatio-temporal chaos appears
\cite{Doering1988,Malomed1990,Tuckerman1990}. We will now use the
result about the amplitude equation in order to trace back this destabilization
transition to the original coupled oscillator system of the form (\ref{eq:main}).
Specific attention will be paid to the scaling behavior of the parameter
region where this transition takes place. 

Starting with the scale free amplitude equation 
\begin{equation}
\partial_{T_{2}}u=r\kappa_{2}u+\frac{\kappa_{3}}{2}\partial_{\xi}^{2}u+\zeta u\left|u\right|^{2}\label{eq:22-1}
\end{equation}
with periodic boundary conditions, we assume that at $r=r_{0}$ a
destabilization of the homogeneous stationary state takes place and
for further increasing parameter $r$ a transition to spatio-temporal
chaos can be observed. This transition takes place within some bounded
interval $\Delta r.$ Considering now the corresponding behavior of
an oscillator system that is described by this amplitude equation,
we have to introduce the scalings given in Proposition \ref{prop2}.
In particular, according to the parameter scaling $p=r\varepsilon^{2}$,
the same transition takes place in the coupled oscillator system within
a parameter interval 
\[
\Delta p(N)=\Delta r\varepsilon^{2}=\Delta r\frac{1}{N^{2}},
\]
of the corresponding parameter $p$ of the oscillator system. Therefore,
one can expect that a generic transition to chaos (or hyperchaos)
in the system of $N$ coupled oscillators occurs within a parameter
interval that scales as $1/N^{2}$. In the following section, we present
an example where the assumptions that we made above are satisfied
and, consequently, in the oscillator system an almost immediate transition
from stationary behavior to chaos can be observed.

\section{Example: a ring of unidirectionally coupled Duffing oscillators}

\begin{figure}[ht]
\begin{centering}
\includegraphics[width=11cm]{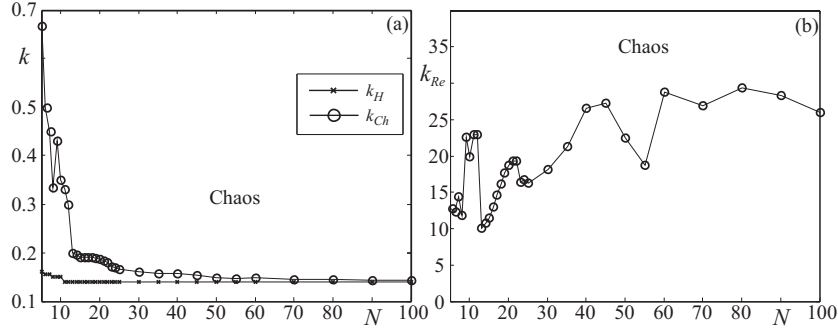} 
\par\end{centering}

\caption{\label{fig:Values-of-coupling}(a) Values of the coupling parameter
$k_{H}$ at Hopf bifurcation (crosses) and $k_{Ch}$ at the transition
to chaos (circles) for increasing number of oscillators $N$. (b)
Rescaled transition interval $k_{Re}$(circles) versus $N$. }
\end{figure}

In the previous section we discussed the scaling law for the transition
to chaos in a ring of coupled oscillators. In this section we provide
a numerical example to illustrate this scaling. We use the autonomous
Duffing oscillator described by the following second order ordinary
differential equation

\begin{equation}
\ddot{y}+d\dot{y}+ay+y^{3}=0,\label{eq:duff_single}
\end{equation}
where $d$ and $a$ are positive constants. System (\ref{eq:duff_single})
is a single-well Duffing oscillator, which has a single equilibrium
point at $y=\dot{y}=0$. Due to the presence of damping $(d>0)$ this
equilibrium is an attractor for all initial conditions. We consider
now a ring of $N$ such oscillators with a linear unidirectional coupling
to the next neighbor. Introducing new coordinates $x=y$, $z=\dot{y}$
and the coupling into Eq. (\ref{eq:duff_single}), the equations of
motion have the form

\begin{equation}
\begin{array}{l}
\dot{x_{j}}=z_{j},\\
\dot{z}_{j}=-dz_{j}-ax_{j}-x_{j}^{3}+k\left(x_{j+1}-x_{j}\right),
\end{array}\label{eq:duff_sprz_schemat}
\end{equation}
 where $k$ is the coupling coefficient and indices are considered
modulo $N$. A detailed analysis of the behavior of this system was
already presented in \cite{Perlikowski2010}. For at least three coupled
oscillators and increasing coupling strength $k$, one can observe
rich dynamics starting from periodic oscillations to hyperchaos. Here
we focus our attention on the transition to chaos (hyperchaos) and
its dependence on the number of oscillators. In our numerics, we used
the fixed parameter values $a=0.1$ and $d=0.3$. We calculated the
maximum Lyapunov exponents for increasing $N$ and varying coupling
coefficient $k$. based on this, we determined two values: $k_{H}$
where the uniform stationary state loses its stability in a Hopf bifurcation,
and $k_{Ch}$ where the transition to chaos takes place. A reliable
computation of the Lyapunov exponents was possible only for moderate
values of $N$ ($N<30$). For larger systems, we monitored the behavior
in an appropriately chosen Poincare section in order to determine
$k_{H}$ and $k_{Ch}$. Fig.~\ref{fig:Values-of-coupling}(a) shows
that for increasing $N$ the distance between $k_{H}$ and $k_{Ch}$
decreases. Finally, for a large enough number of oscillators, the
transition to chaos appears practically immediately after the Hopf
bifurcation. In order to validate our conjecture about the $1/N^{2}$
scaling, we plot in Fig.~\ref{fig:Values-of-coupling}(b) the scaled
transition intervals

\[
k_{Re}=\left(k_{Ch}-k_{H}\right)N^{2}.
\]
The numerical results indicate that the scaled transition interval
neither tends to zero nor diverges to infinity, which supports the
scaling assumption.

\begin{figure}[ht]
\begin{centering}
\includegraphics[width=11cm]{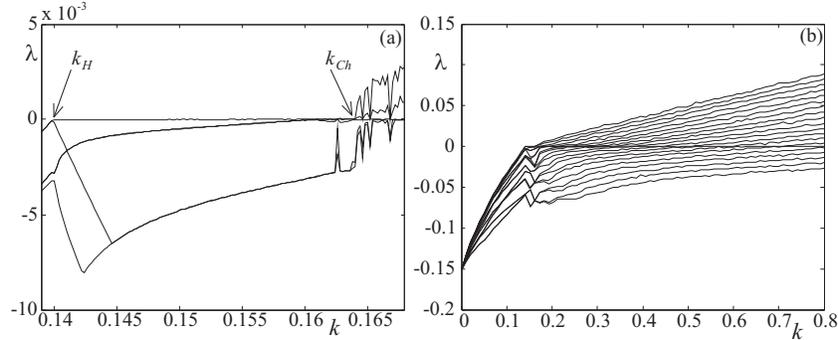} 
\par\end{centering}

\caption{\label{fig:lap_sec30}Lyapunov exponents for $N=30$ unidirectionally
coupled Duffing oscillators (a) five largest Lyapunov exponents, $k_{H}$
and $k_{Ch}$ indicate the Hopf bifurcation and the transition to
chaos, $k\in(0.139,0.17)$, (b) twenty largest Lyapunov exponents
for $k\in(0.0,0.8)$. At $k=0.8$, there are fourteen positive Lyapunov
exponents. }
\end{figure}

In Fig. \ref{fig:lap_sec30} we present the Lyapunov spectrum for
thirty coupled oscillators. Panel (a) shows that at $k_{Ch}=0.164$
two Lyapunov exponents become positive practically simultaneously.
A similar behavior has been observed for other large values of $N$.
This shows that the transition to chaos and to hyperchaos occursin
a large ring almost simultaneously. With further increasing coupling
coefficient $k$ (Fig. \ref{fig:lap_sec30}(b)) one can observe that
more and more Lyapunov exponents become positive which leads to a
high dimensional chaos \cite{Kapitaniak1991,Kapitaniak1993a}.

\section{Conclusions}

The results of this paper are twofold. The first one is more theoretical:
we have derived an amplitude equation for small amplitude solutions
of a ring system of coupled oscillators. This equation has the form
of a complex Ginzburg-Landau equation with periodic boundary conditions.
Such an amplitude equation is typical for reaction-diffusion systems
\cite{Eckhaus1965,Tuckerman1991} in continuous media and can also
be found for systems with large delay \cite{Wolfrum2006}. Obviously,
a similar behavior should be expected in discrete systems with diffusive
coupling. However, our results show that in coupled oscillator systems
the Ginzburg-Landau equation can also be used to describe the dynamics
in systems with unidirectional, i.e. purely convective coupling or
even with anti-diffusive interaction. In this sense, the class of
coupled oscillators that we treated in this paper differs substantially
from discrete analogs of the classical results for continuous media.

For the second result, we applied the amplitude equation to the specific
scenario of the transition to spatio-temporal chaos in such systems.
In this way, we provided a general framework for the coupling induced
transition to high dimensional chaos in coupled oscillator systems,
that has been described for a specific example in \cite{Perlikowski2010}.
In particular, we were able to show that the observed fast transition
to chaos is a generic feature and follows a $1/N^{2}$ scaling law
for systems with a large number $N$ of oscillators. As a result of
this scaling, one observes for large systems a practically immediate
transition from a uniform stationary state to chaos. We illustrated
this behavior by a numerical example, where the emergence of high
dimensional chaos is demonstrated by corresponding Lyapunov spectra.

\begin{acknowledgments}
S. Yanchuk acknowledges the support of the DFG collaborative research
center ''Control of self-organizing nonlinear systems: Theoretical
methods and concepts of application'' (SFB910) under the project A3.
P. Perlikowski, A. Stefa\'{n}ski and T. Kapitaniak acknowledge the
support from Foundation for Polish Science - Team Programme (Project
No TEAM/2010--5/5) and P. Perlikowski START Programme. 
\end{acknowledgments}

\bibliographystyle{aipnum4-1}
%
\end{document}